\documentclass[11pt]{amsart} 
\usepackage{amsfonts,amscd}
\def\Z{{\mathbb Z}}

\newtheorem{thr}{Theorem}[section]
\newtheorem{df}[thr]{Definition}
\newtheorem{lm}[thr]{Lemma}
\newtheorem{pr}[thr]{Proposition}

\newtheorem{thr*}{Theorem}
\newtheorem{pr*}{Proposition}
\newtheorem{co*}{Corollary} 
\begin{document} 
\baselineskip 22pt
\title[Smooth structure of some symplectic surfaces]{Smooth structure of some symplectic surfaces}
\author{Stefano Vidussi}
\address{Department of Mathematics, University of California, Irvine, California 92697 }
\email{svidussi@math.uci.edu}
\maketitle
\baselineskip 18pt
\section{Introduction} McMullen and Taubes \cite{McMT} have constructed a
remarkable simply connected smooth $4$-manifold, denoted by $X$,  
starting from a $4$-component link $K \subset S^{3}$ and
four copies of the rational elliptic surface $E(1)$. The interest in the link $K$ stems from the fact that
it admits several inequivalent fibrations over $S^{1}$; these inequivalent fibrations give rise to two
inequivalent symplectic structures on $X$, providing the first simply connected example of manifold with
this property. The ingredients in the construction of \cite{McMT} are reminiscent of those used by Fintushel
and Stern in defining a large class of smooth $4$-manifolds, and it is natural to ask how these
constructions are related. In this note we will compare the link surgery construction of \cite{FS} and the
McMullen-Taubes example in order to prove that the latter manifold is 
diffeomorphic to a Fintushel-Stern
manifold. This analysis (further developed in \cite{V}) will lead us to introduce a new presentation of $X$
that allows us to identify  a new symplectic structure on $X$.
We will assume some familiarity with \cite{FS} and \cite{McMT}.
\section{Construction of the 4-manifolds} \label{prime}
We start by recalling the link surgery construction of \cite{FS}, omitting
(for the sake of brevity) full
generality. Consider an $n$-component oriented link $K \subset S^{3}$. Let 
$p_{i} = - \sum_{j \neq i} lk(K_{i},K_{j})$; the closed manifold $M_{K}$ obtained by performing
$p_{i}$-surgery on the $i$th component has the property that the image $m_{i}$ of
each meridian $\mu(K_{i})$ has infinite order in $H_{1}(M_{K},\Z)$ and is
canonically framed; in $S^{1} \times M_{K}$, the tori $S^{1} \times
m_{i}$ have self-intersection zero and are framed and essential
in homology. Next take $n$ copies of the simply connected elliptic surface without multiple fibers
$E(m)$, each containing an elliptic fiber $F_{i}$, and
construct, by normal connected sum, the manifold 
\begin{equation} \label{fs} E(m)_{K} = \coprod E(m)_{i} \#_{F_{i} =
S^{1} \times m_{i}} S^{1} \times M_{K}. \end{equation}
The gluing is made so as to send the homology class of the
  normal circle to the $i$th torus $S^{1} \times m_{i}$, represented by $p_{i}m_{i} + l_{i}$ (where
$l_{i}$ is the image of the preferred longitude $\lambda(K_{i})$) to the class of a normal circle to the
$i$th elliptic fiber. These prescriptions can be insufficient to uniquely define the manifold:
the gluing map is defined up the action of $SL(3,\Z)$ matrices of the form 
\begin{equation}
\left( \begin{array}{ccc} a & b & 0 \\ d & e & 0 \\ g & h & 1
  \end{array} \right); \end{equation} since $F$ is in
the neighborhood of a
cusp fiber in $E(m)$, we can dispose of the indeterminacy
corresponding to the upper left $SL(2,\Z)$ factor (due to the
absence of a canonical choice for the basis of $H_{1}(F,\Z)$) because any
fiber and orientation preserving diffeomorphism of $\partial (E(m)
\setminus \nu F)$ extends to a (fiber-preserving) diffeomorphism of
$E(m) \setminus \nu F$ (see Chapter 8 of \cite{GS}); the symbol $\nu (\cdot)$ denotes the open 
neighborhood of an embedded submanifold. 
The remaining indeterminacy, however, cannot be disposed of in general.
The manifold $E(m)_{K}$ is simply connected and has $b^{+}_{2} \geq n$.
\par We will discuss now the example of McMullen and Taubes.
 Consider, in $S^{3}$, the 4-component oriented link $K$
given by the
union of the Borromean rings $K_{1} \cup K_{2} \cup K_{3}$ 
and the axis of $\Z_{3}$-symmetry $K_{4}$. Let
$N := S^{3} \setminus \nu K$. We recall the form of the Alexander polynomial 
$\Delta_{K}(x,y,z,t)$ of $K$; here
$x,y,z$ are the variables corresponding 
to the meridians of the Borromean rings and $t$ corresponds to the
meridian to the axis:
\begin{equation} \label{poly} \begin{array}{c} \Delta_{K}(x,y,z,t) = -4 + (t +
    t^{-1}) + (x + x^{-1} + y + y^{-1} + z + z^{-1}) + \\ \\ - (xy +
  (xy)^{-1} + yz + (yz)^{-1} + xz + (xz)^{-1}) + (xyz + (xyz)^{-1}). 
\end{array} \end{equation} We have another description for
$N$; perform $0$-surgery on $S^{3}$ along the components of the 
Borromean rings; it is well known that this surgery yields $T^{3}$.
We can thus write $N = S^{3} \setminus \nu K = T^{3} \setminus \nu L$,
where $L$ is a framed link in $T^{3}$, whose first three components give a basis of $H_{1}(T^{3},\Z)$.
In fact, when we perform the 
$0$-surgery on the Borromean rings, the three
meridians $\mu(K_{i})$ ($i = 1,2,3$) to the components
of the Borromean rings go over longitudes $m_{i}$ of $L_{i}$, while
the preferred longitudes $\lambda(K_{i})$ are sent to 
meridians $l_{i}$ of $L_{i}$. The longitude of
$K_{4}$ becomes a longitude
to the component $L_{4} \subset T^{3}$, which satisfies the
relation $L_{4} = L_{1} + L_{2} + L_{3}
\in H_{1}(T^{3},\Z)$; the meridian $\mu(K_{4})$ of $K_{4}$ goes instead
to a meridian $m_{4}$ of $L_{4}$ and is null-homotopic in $T^{3}$.
It is instead nontrivial in $H_{1}(N,\Z)$, where the four generators
are given by the meridians. We have $H^{1}(N,\Z) \supset i^{*}H^{1}(T^{3},\Z) = \Z<t>^{\perp}$.
Then define the normal
  connected sum \begin{equation} \label{xman} X = \coprod E(1)_{i} \#_{F_{i} =
    S^{1} \times L_{i}} S^{1} \times T^{3}. \end{equation}
Again, the definition requires that the homology class of the normal circle 
to $S^{1} \times L_{i}$ be
sent to the homology class of the normal circle to the $i$th elliptic fiber. 
The previous remarks on the ambiguity of the definition apply. This
manifold is simply connected and has $b^{+}_{2} > 1$.
\par We show now that both constructions appear as
particular cases of a general construction: consider the exterior of an oriented $n$-component link $K
\subset S^{3}$ together with the choice, in each boundary component, of an homology basis of simple curves
$(\alpha_{i},\beta_{i})$ of intersection
$1$. We introduce the following definition.
\begin{df} \label{emk} Take a link $K$ as above with homology basis $(\alpha_{i},\beta_{i})$ and an
elliptic surface
$E(m)$. Define the manifold \begin{equation} \label{descr}
E(m;\alpha_{i},\beta_{i}) = (\coprod E(m)_{i} \setminus \nu F_{i}) \cup_{F_{i} \times S^{1} = S^{1}
\times
\alpha_{i} \times \beta_{i}} (S^{1} \times (S^{3} \setminus \nu K)) \end{equation} where
the gluing is made by lifting a diffeomorphism between $S^{1} \times
\alpha_{i}$ and $F_{i}$ to an orientation-reversing
diffeomorphism of the boundary tori in such a way
that the homology class of $\beta_{i}$ is sent to the homology class
of the normal circle to the $i$-th
  elliptic fiber. \end{df} \noindent The gluing condition is not enough to 
define the manifold completely. 
As in the case of Fintushel-Stern manifolds, the
  ambiguity related to the absence of a chosen
  basis in $H_{1}(F_{i},\Z)$ is only apparent, whereas the remaining
  ambiguity is effective. Moreover, the smooth manifold (as the notation
  suggests) can depend on the choice of the $(\alpha_{i}, \beta_{i})$,
  with the noteworthy exception considered in the following lemma. \begin{lm}
    \label{coin} Let
  $E(1;\alpha_{i},\beta_{i})$ be defined as before.
  Then the manifold is well-defined and moreover its definition
  depends uniquely on $K$: that is, it is unaffected by the choice of the
  basis on $\partial (S^{3} \setminus \nu K)$. \end{lm} \noindent {\bf Proof:} 
 This follows from the fact that {\it any} orientation-preserving 
diffeomorphism of $\partial
(E(1) \setminus \nu F)$, and not only the fiber-preserving ones, extends to 
an orientation-preserving
diffeomorphisms of $(E(1) \setminus \nu F)$ (see \cite{GS}):
on each boundary
  component we can reabsorb any orientation-preserving 
self-diffeomorphism of $S^{1} \times
  \alpha_{i} \times \beta_{i}$ by an orientation-preserving self-diffeomorphism
of $\partial
(E(1)_{i} \setminus \nu F_{i})$, which extends to $E(1)_{i}
  \setminus \nu F_{i}$.
No matter how we glue the
manifold $S^{1} \times (S^{3} \setminus \nu K)$ (in particular, for any choice of homology basis for the
boundary), the resulting four manifolds are smoothly equivalent. \endproof 
\noindent Analyzing the previous construction yields the following, 
straightforward proposition. \begin{pr} The
Fintushel-Stern manifolds
$E(m)_{K}$ and the McMullen-Taubes manifold $X$ can be described via the construction in Definition \ref{emk}.
\end{pr} 
\noindent {\bf Proof:} The definition of normal connected sum
shows that the manifolds defined in equation \ref{fs}
can be rewritten in the form 
\begin{equation} \label{oldglue}
E(m)_{K} = ( \coprod E(m)_{i} \setminus \nu F_{i}) \cup
(S^{1} \times (S^{3} \setminus \nu K)) \end{equation}
where the gluing is made lifting a diffeomorphism between $S^{1} \times
\mu(K_{i})$ and $F_{i}$ to an orientation-reversing
diffeomorphism of the boundary tori so that the homology class of $p_{i}
\mu(K_{i}) + \lambda(K_{i})$ is sent to the class of the normal circle
to $F_{i}$. Hence the manifold $E(m)_{K}$ corresponds to the choice 
  $(\alpha_{i},\beta_{i}) = (\mu(K_{i}),p_{i} \mu(K_{i}) +
  \lambda(K_{i}))$. 
 \\ Concerning the McMullen-Taubes example, 
an analysis of the definitions via normal
  connected sum of Eq. \ref{xman} 
(keeping track of the framing of $L_{i}$) 
shows, as $T^{3} \setminus \nu L = S^{3}
\setminus \nu K$,
that $X$ corresponds to $m=1$ and to the choice 
$(\alpha_{i}, \beta_{i}) = (\mu(K_{i}),\lambda(K_{i}))$ for $i = 1,2,3$ and
$(\alpha_{4},\beta_{4}) = (\lambda(K_{4}),-\mu(K_{4}))$. 
\endproof
\noindent Note that the latter definition
differs from the Fintushel-Stern one, applied to the same link, 
for the different choice of the homology
basis. However, in this particular case, we have our next lemma.
\begin{lm} \label{mine} The McMullen-Taubes manifold
  $X$ is diffeomorphic to the Fintushel-Stern manifold
  $E(1)_{K}$. \end{lm} \noindent {\bf Proof:} This follows as particular case of
  Lemma \ref{coin}. The same argument implies also that the manifold is well defined. \endproof 
\section{symplectic structures}
We now want compare the symplectic structure arising naturally from the different presentations of $X$. 
The proof of the existence of symplectic structures on $X$ follows by application of Gompf's theorem on
symplectic normal connected sum between $\coprod_{i} E(1)$ and 
$S^{1} \times M_{K}$ (resp., $S^{1} \times T^{3}$)  in the
Fintushel-Stern (resp., McMullen-Taubes) construction.
Both $M_{K}$ and $T^{3}$ are fibered $3$-manifolds obtained by Dehn filling 
of $S^{3}
\setminus \nu K$ along the different surgery curves. For any choice of a 
fiber $\Sigma$ in $M_{K}$ (resp.,
$T^{3}$) transverse to the image of the link, $E(1)_{K}$ (resp.,
$X$) inherits a natural symplectic structure induced from the closed, nondegenerate $1$-form
defining the fibration on $S^{3} \setminus \nu K$. For any link $K$, fibrations on $S^{3} \setminus \nu
K$ are identified with the elements of $H^{1}(S^{3} \setminus \nu K,\Z)$ laying on the cones over some of
the top dimensional faces of the Thurston unit sphere. The latter is defined, for $\varphi \in H^{1}(S^{3}
\setminus
\nu K,\Z)$, by minimizing the quantity
\begin{equation} \chi(\Sigma) = \sum_{\chi(\Sigma_{i}) < 0}(- \chi(\Sigma_{i})) \end{equation} among
properly embedded representatives $\Sigma$ of the Poincar\'e dual of 
$\varphi$ and then extending by linearity and continuity to real 
cohomology classes.     
The fibration on $M_{K}$ restricts by construction (see \cite{FS}) to the fibration of
$S^{3} \setminus \nu  K$ with fiber given by the minimal spanning surface of the link $K$, that is, to the
class $(1,1,1,1)
\in H_{2}(S^{3} \setminus \nu K,\Z)$. On $T^{3}$, as discussed in \cite{McMT}, every fibration that
restricts to the cone over the top-dimensional faces of the Thurston unit sphere on $i^{*}H^{1}(T^{3},\Z)
\subset H^{1}(S^{3} \setminus \nu K,\Z)$ induces a symplectic structure on $X$. We can relate
the fibration of class $(1,1,1,1)$ and the fibrations laying in 
$i^{*}H^{1}(T^{3},\Z)$: the analysis of the
Thurston norm on $H^{1}(S^{3} \setminus \nu K,\Z)$, detailed in \cite{McMT}, shows that $(1,1,1,1)$ lies in
the cone over the top dimensional face identified by the dual vertex $xyz$ (we use the same notation as
Eq. \ref{poly}), a face that already contains fibered elements of 
$i^{*}H^{1}(T^{3},\Z)$. As a consequence, the
canonical bundle corresponding to the symplectic structure induced on $X$ by this fibration cannot be used to
distinguish it from the ones exhibited in \cite{McMT}. 
\\ Let's discuss now how we can produce a new symplectic structure that can be
distinguished from the known ones by studying the canonical class. 
The unit sphere of the Thurston norm of $S^{3} \setminus \nu K$ is given, as discussed in
\cite{McMT}, by the product of the unit sphere in the subspace $i^{*}H^{1}(T^{3},\Z)$ and the interval
$[-\frac{1}{2},\frac{1}{2}]$ of the orthogonal subspace: every fibered face is determined by a dual
vertex among the sixteen vertices of the Newton polyhedron of the Alexander 
polynomial.
We can represent the orthogonal subspace to $i^{*}H^{1}(T^{3},\Z)$ as pull
back under inclusion of the first cohomology group of $S^{1} \times S^{2}$: in fact, $0$-surgery on the axis
$K_{4}$ of the Borromean ring exhibits $N$ as complement of a link ${\hat L}$ in
$S^{1} \times S^{2}$. The images of the meridians $\mu(K_{i})$ for $i=1,2,3$ are (nullhomotopic) meridians 
to the components of ${\hat L}$ with the same index; $\mu(K_{4})$ goes to a preferred longitude of ${\hat
L}_{4}$. The longitudes $\lambda(K_{i})$ for $i=1,2,3$ go to preferred
longitudes of the respective ${\hat L}_{i}$, while $\lambda(K_{4})$ goes to a meridian to ${\hat L}_{4}$. 
The fiber of $S^{1} \times S^{2}$ restricts to
the fiber of $S^{3} \setminus \nu K$ identified by the homology class $(0,0,0,1) \in H^{1}(S^{3}
\setminus \nu K,\Z)$ (a disk spanning the axis, pierced once by each components of the Borromean rings). 
We have now the following. \begin{df}
\label{defy} Consider the framed symplectic tori $S^{1} \times {\hat L}_{i} \subset S^{1} \times S^{1}
\times S^{2}$  of  self-intersection zero together with 
four copies of the rational elliptic surface $E(1)$. We define
the normal connected sum 
\begin{equation} \label{myma} Y = \coprod E(1)_{i} \#_{F_{i} =
    S^{1} \times {\hat L}_{i}} S^{1} \times S^{1} \times S^{2}. \end{equation}
The definition of normal connected sum imposes
that the homology class of the normal circle to $S^{1} \times {\hat L}_{i}$ be
sent over the homology class of the normal circle to the $i$-th elliptic fiber. \end{df}
\noindent This definition yields immediately our next proposition. 
\begin{pr} The manifold $Y$ introduced in Definition
\ref{defy} is a manifold of type $E(1;\alpha_{i},\beta_{i})$; it is,
moreover, diffeomorphic to $X$ and to the
Fintushel-Stern manifold $E(1)_{K}$. \end{pr} \noindent {\bf Proof:} The first statement follows by observing
that the definition of $Y$ corresponds  to the choice $S^{1} \times S^{2} \setminus \nu {\hat L} = S^{3}
\setminus
\nu K$ and to the homology basis $(\alpha_{i},\beta_{i}) = (\lambda(K_{i}),-\mu(K_{i}))$ for $i=1,2,3$ and
$(\alpha_{4},\beta_{4}) = (\mu(K_{4}),\lambda(K_{4}))$. The second statement is a corollary, as Lemma
\ref{mine}, of Proposition \ref{coin}.
\endproof \noindent The construction of $X$ introduced in Def. \ref{defy} induces naturally a
symplectic structure on the manifold: the fibration of $S^{3} \setminus \nu K$ with class $(0,0,0,1)$ has dual
vertex
$t$, as we can see by looking at the Alexander polynomial in equation \ref{poly}. Theorem 3.4 of
\cite{McMT} identifies the canonical bundle of this symplectic structure as the image
of twice this vertex under the injective map $H_{1}(S^{3} \setminus \nu K,\Z) \rightarrow H^{2}(X,\Z)$. This
canonical bundle has different valence, as vertex of the Newton polyhedron of the SW polynomial, than the
canonical bundles obtained from the previous two construction of $X$ (see \cite{McMT}) and so is
combinatorially different. As a consequence, it lies in a different orbit with repect to
the action of the diffeomorphism group of $X$, that acts by preserving the Newton polyhedron.   
This proves the following. \begin{thr} The symplectic structure induced by
normal connected sum on $Y$ is not equivalent (up to combination of pull-back and homotopies) to the previous
ones.
\end{thr}
\noindent The Seiberg-Witten polynomial of $X$ is given by $\Delta_{K}(x^{2},y^{2},z^{2},t^{2})$; the new
symplectic structure (and its conjugate), together with the fourteen constructed in \cite{McMT}, exhaust the sixteen
basic classes with coefficient $\pm 1$.
\par In \cite{V} we discuss how the above constructions can be extended to
obtain further generalizations of the Fintushel-Stern link surgery 
construction.

\end{document}